\documentclass[a4paper,12pt]{amsart}\usepackage{amsfonts,amssymb,oldgerm,enumerate,mathrsfs,bbm,hyperref,rotating}
\UseRawInputEncoding

\def\C{\mathbb{C}}\def\bbf{{\bar{f}}}

\def\P{\mathbb{P}}\def\R{\mathbb{R}}\def\bX{{\bar{X}}}

\def\di{\partial}

\renewcommand{\stackrel}[2]{\ \lower 0.2ex \hbox{$\mathrel{\mathop{#2}\limits^{#1}}$}\ }

\def\De{\Delta}

\newcommand{\bbm}{\begin{bmatrix}}\newcommand{\ebm}{\end{bmatrix}}
\newcommand{\ber}{\begin{array}{l}}\newcommand{\eer}{\end{array}}
\newcommand{\bpm}{\begin{pmatrix}}\newcommand{\epm}{\end{pmatrix}}
\newcommand{\bM}{\begin{matrix}}\newcommand{\eM}{\end{matrix}}
\newcommand{\bee}{\begin{enumerate}}\newcommand{\eee}{\end{enumerate}}
\newcommand{\bei}{\begin{itemize}}\newcommand{\eei}{\end{itemize}}

\def\sset{\subset}\def\sseteq{\subseteq}\def\smin{\setminus}

\newcommand{\beq}{\begin{equation}}\newcommand{\eeq}{\end{equation}}
\newtheorem{Lemma}{Lemma}[section]\newcommand{\bel}{\begin{Lemma}}\newcommand{\eel}{\end{Lemma}}
\newtheorem{Example}[Lemma]{Example}\newcommand{\bex}{\begin{Example}\rm}
\newcommand{\eex}{\end{Example}}
\newtheorem{Proposition}[Lemma]{Proposition}\newcommand{\bprop}{\begin{Proposition}}\newcommand{\eprop}{\end{Proposition}}
\newtheorem{Property}[Lemma]{Property}\newcommand{\bproperty}{\begin{Property}}\newcommand{\eproperty}{\end{Property}}
\newtheorem{Definition-Proposition}[Lemma]{Definition-Proposition}

\def\bpr{~\\{\em Proof.\ }}
\newcommand{\epr}{{\hfill\ensuremath\blacksquare}\\}
\newtheorem{Theorem}[Lemma]{Theorem}\newcommand{\bthe}{\begin{Theorem}}\newcommand{\ethe}{\end{Theorem}}
\newtheorem{Definition}[Lemma]{Definition}\newcommand{\bed}{\begin{Definition}}\newcommand{\eed}{\end{Definition}}
\newtheorem{Remark}[Lemma]{Remark}\newcommand{\beR}{\begin{Remark}\rm}\newcommand{\eeR}{\end{Remark}}
\newtheorem{Corollary}[Lemma]{Corollary}\newcommand{\bcor}{\begin{Corollary}}\newcommand{\ecor}{\end{Corollary}}
\newcommand{\bet}{\begin{tabular}{cccccccc}}\newcommand{\eet}{\end{tabular}}

\def\into{\stackrel{j}{\hookrightarrow}}

\title[]{A\MakeLowercase{ndreotti-}F\MakeLowercase{rankel-}H\MakeLowercase{amm theorem for morphisms of algebraic varieties}}
\author[]{D\MakeLowercase{mitry} K\MakeLowercase{erner}}
\address{Department of Mathematics, Ben Gurion University of the Negev, P.O.B. 653, Be'er Sheva 84105, Israel. dmitry.kerner@gmail.com}

\subjclass[2020]{Primary
14F45	
\quad Secondary
14B05	
32C18	
32S50	
}

\keywords{Topology of complex affine varieties, Topology of Stein spaces, Singularity Theory, Topology of Milnor Fibre}

\date{\today\ \  filename: \jobname.tex}
\thanks{I was supported by the Israel Science Foundation, grant No.  1405/22}
\begin{document}
 \maketitle
  \begin{abstract}
The classical Andreotti-Frankel-Hamm theorem reads: a complex affine algebraic variety  $B,$ of $dim_\C B=n,$ has homotopy type of $dim_\R\le n.$
 We prove the relative version   for morphisms $X\to B.$
 \end{abstract}

\section{Introduction}
\subsection{}
Let   $B$ be a complex affine algebraic variety  or a Stein space (with arbitrary singularities),   of $dim_\C B=n.$ The classical theorem reads:
\beq\label{Eq.AFKHM}
 \text{\em $B$  admits a  deformation-retraction to a closed subanalytic subspace    of $dim_\R\le n.$}\vspace{-0.4cm}
 \eeq
\[\text{\em Moreover, if $B$ is affine algebraic, then the resulting subspace can be chosen compact.}
\] 
Apparently the first result of this type was   in  \cite{Andreotti-Frankel}: a Stein manifold
 $B$ has the homotopy type of a CW-complex  of real dimension $dim_\R\le n.$
     The statement was strengthened, in Theorem 7.2 of \cite{Milnor}, to the deformation-retractions of Stein manifolds.
     The statement on the homotopy type was extended to simply connected affine algebraic sets (with arbitrary singularities) in \cite[pg.49]{Kato}, then to arbitrary  affine algebraic sets in Theorem 2.12 of \cite{Karchyauskas}, and finally to  Stein spaces
 (with arbitrary singularities) in  \cite{Hamm.83}. The version  \eqref{Eq.AFKHM} of deformation retraction   is Theorem 1.1 in  \cite{Hamm-Mihalache}.

\medskip

This theorem is an everyday tool for  Geometry and Topology of complex varieties and Stein spaces. In Singularity Theory it implies: the Milnor fibre of any singular germ has homotopy type of $dim_\R\le n.$ This fundamental property is the starting point in the study of Milnor fibres, a flourishing field in the last 70 years, see e.g. Chapters 6,8,9  in \cite{Handbook.I} and Chapters 6,7 in \cite{Handbook.II}.

\subsection{}
   Take a dominant morphism of complex algebraic varieties $X\stackrel{f}{\to} B,$ think of it as a fibration over the base $B.$ This brings two questions of A.F.K.H.-type:
   \bei
   \item (vertical version) Suppose all the fibres of $X\to B$ are affine algebraic varieties of $dim_\C\le n.$ Does the family admit a fibrewise deformation-retraction to a subfamily $Z\to B,$ with all fibres of $dim_\R\le n?$

   \item (horizontal version) Suppose    $B$ is affine, $dim_\C B\le n.$ Does there exist a deformation retraction $B\rightsquigarrow B',$ with $dim_\R B'\le n,$ that lifts to a deformation-retraction of the whole family, sending fibres to fibres?
   \eei
   (Of course, in both cases the deformation-retraction is non-analytic.)

   The vertical version is trivial. Take any stratification making $f$ a $C^0$-stratified-trivial morphism, and invoke A.F.K.H. successively over the strata. For a significantly stronger statement see \cite{Hamm.01}.

   The horizontal version is less trivial. Not many deformation-retractions of $B$ can  be lifted to  deformations of $X,$ because $f$ is not a $C^0$-locally trivial fibration.

\medskip 

 We construct a special deformation-retraction of $B,$ establishing the horizontal version.
 \\
 \parbox{12cm}{
 \bthe\label{Thm.A.F.K.H.for.morphisms}
  There exists a deformation-retraction $\{\Phi^B_t\}$ that
 lifts to $X,$ making the commutative diagram. The lifted version satisfies:  \hspace{0.8cm} $\Phi_0=Id_X,\hspace{0.8cm} dim_\R  \Phi^B_1(B)\le n,\hspace{0.8cm} \Phi_t|_{\Phi_1(X)}=Id.
 $
 \ethe
}\hspace{1cm}
$\bM
X\times[0,1]\stackrel{\{\Phi_t\}}{\to}X
\\(f,Id)\downarrow \quad\quad\quad \downarrow f\\
B\times[0,1]\stackrel{\{\Phi^B_t\}}{\to}B
\eM$

This homotopy is (naturally) called: a deformation-retraction of the morphism $f.$

While the statement of Theorem \ref{Thm.A.F.K.H.for.morphisms} is quite natural, it appears  to be not known.

\

Our proof is an iterative argument on bifurcation loci, and uses the A.F.K.H. theorem (for algebraic varieties) as a ``black-box".
 E.g. in the trivial case, $X=B$ and $f=Id_B,$ the proof just refers to the original statement(s).
 We do not address the case of Stein morphisms, as their bifurcation loci can be quite pathological.

 \

This relative version of the classical A.F.K.H. theorem is of immediate use e.g. in Algebraic Geometry, in the study of families and fibrations.
 In Singularity Theory it is needed e.g.  to study the Milnor fibre via its projections (or via other maps).

 Our particular motivation was the study of fast vanishing cycles on links of singular germs, \cite{Kerner-Mendes.Coverings.Discriminants},
  \cite{Kerner-Mendes}.

\subsection{Acknowledgements} The results were obtained mainly during the conferences ``Metric Theory of Singularities" (Krakow, October 2024), and ``Logarithmic and non-archimedean methods in singularity theory" (Luminy, January 2025). Thanks to their organizers.

\section{Several general (topological) facts}\label{Sec.General.Facts} All the spaces below are finite CW complexes. By $dim_\R X$ we mean the maximum of the dimensions of cells. A deformation-retraction to a subset, $X\rightsquigarrow X',$ is a homotopy of $X$ to $X'$ that restricts to identity on $X'.$

 A space homotopically equivalent to a space of $dim_\R\le n$ is called ``of homotopy type $dim_\R\le n$".
\bee[\bf i.]
\item Take a pair $X\supset Y$ of contractible spaces with the homotopy extension property. Then $X$ admits a deformation-retraction to $Y.$
  See e.g. Corollary 0.20 of \cite{Hatcher}.

 All our pairs will be CW pairs, thus having the homotopy extension property, by Proposition 0.16 of \cite{Hatcher}.

 \item Suppose $X$ is of homotopy type of $dim_\R\le n.$ Then $X\cup Cells_{n+1}$ is contractible, where $ Cells_{n+1}$ denotes a union of cells of $dim_\R\le n+1,$ and they are glued to $X$ along their boundaries.  This follows, e.g. by the Whitehead theorem, \cite[pg.346]{Hatcher}.

\item
Suppose a space $Y\cup Z$ is contractible, and its subspaces $Z,W$ are contractible. Suppose $dim_\R(Y\cap Z\cap W)< dim_\R Z -1.$ Then there exists a deformation retraction $Y\cup Z\rightsquigarrow Z\cup W\cup \De,$ where $\De\sset \overline{Y\smin(Z\cup W)}$ and  $dim_\R\De< dim_\R Z.$

\bpr
If $Z\cup W$ is already contractible, then we invoke \S\ref{Sec.General.Facts}.i. for the pair
 $Z\cup W\sset Z\cup Y.$ In the general case $Z\cup W$ has a non-trivial homology and fundamental group.
 We eliminate these   cycles by gluing cells. Once all the cycles are cancelled, we define $\De$ as the union of these cells, and  invoke \S\ref{Sec.General.Facts}.ii. Here are the details.

 Recall the Mayer-Vietoris sequence, $H_i(Z)\oplus H_i(W)\to H_i(Z\cup W)\to H_{i-1}(Z\cap W).$ In our case it implies: the non-trivial homology of $Z\cup W$ comes only from that of $Z\cap W.$ If a cycle of $Z\cap W$ has a representative inside $Z\smin Y,$ then it is eliminated by gluing a cell inside $Z.$ (Because $Z$ is contractible.) Such a cell will not contribute to $\De.$ Therefore it is enough to consider those cycles of $Z\cap W$ that lie inside $Z\cap Y.$ But $dim_\R(Z\cap Y\cap W)<dim_\R Z-1,$ thus the cells to be added  are all of dimensions $dim_\R<dim_\R Z.$
  Therefore $dim_\R\De<dim_\R Z.$
  \epr

\item (Thom's first isotopy lemma, e.g. Theorem B.2   in \cite{Mond-Nuno}) Take a stratified topological space, $X=\amalg X_\bullet,$ with $C^1$-smooth strata. Suppose a $C^1$-map $X\stackrel{f}{\to}B$ is proper, and its restriction to each stratum, $f|_{X_j},$ is a submersion (onto the image of the stratum). Then each restriction $f|_{X_j}$ is a $C^0$-locally-trivial fibration.

In our case  $X\stackrel{f}{\to}B$ is a morphism of algebraic varieties, and we make it proper. It admits the trivializing stratification (as above) whose strata $X_\bullet$ are also algebraic varieties.

\item Let $X\stackrel{f}{\to}B$ be a $C^0$-locally-trivial fibration. Then any deformation retraction of $B$ lifts to a deformation-retraction of $f,$ as in the diagram of Theorem \ref{Thm.A.F.K.H.for.morphisms}.
\eee

\section{The proof of the theorem}
\bee[\bf i.]

\item\ [A warmup/the simplest case: $B$ is contractible and $dim_\C B=1.$]

Rename  $B$ to $B_1.$ Take the bifurcation locus, $B_0\sset B_1,$ so that $f$ is  a $C^0$-locally trivial fibration over $B_1\smin B_0.$
  Here $B_0$ is a finite set of points.
 Denote by $Cells_1$ a simply connected path in $B_1,$ connecting  all the points of $B_0.$ Thus $B_0\cup Cells_1$ is contractible inside $B_1.$
  Thus $B_1$ admits a deformation-retraction onto  $B_0\cup Cells_1,$ by \S\ref{Sec.General.Facts}.i.

   This deformation-retraction lifts to a deformation-retraction of the morphism $f:X\to B,$ by the local triviality of $f$ outside of  $B_0\cup Cells_1,$ see \S\ref{Sec.General.Facts}.v. Hence the statement.

\

To extend this argument to higher dimensions we should cope with the locus of local non-triviality of the fibration $X\stackrel{f}{\to}B.$
 As the map $f$ is non-proper, this locus contains  (besides the discriminant of $f$) the set of bifurcations at infinity. This latter can be pathological. Therefore below we compactify the map $f,$ and then work with the discriminant of the proper map.

A deformation retraction of $B$ onto the discriminant (which is identity on the discriminant) lifts to a deformation retraction of $f.$ Then we iterate the argument.

\item\ [Compactification of the map $X\stackrel{f}{\to}B$  to a proper map $\bX\stackrel{\bbf}{\to}B$]

 Embed the source into a projective space, $X\sset \P^N.$ Then re-embed it as the graph $X\into \P^N \times B,$ by $x\to (x,f(x)).$ Thus $j(X)$ is an algebraic subvariety, and we take the closure $\bX:=\overline{j(X)}\sset \P^N \times B.$

 The initial map $f$ acts now by the projection $j(X)\to B.$ Thus it extends to the projection $\bX\stackrel{\bbf}{\to}B.$ This is a proper algebraic morphism.

\item\ [The   trivializing stratification of $\bbf$] Stratify the map, $\bX=\amalg (\bX)_\bullet\stackrel{\bbf}{\to}B.$ Here each stratum $(\bX)_\bullet$ is a smooth algebraic variety, and each restriction $(\bX)_\bullet\stackrel{f|}{\to}B$ is a $C^0$-locally trivial fibration over its image. Moreover, we can take this stratification compatibly with the ``hyperplane at infinity", i.e. with the  splitting $\bX=j(X)\amalg(\bX\smin j(X)) .$ Namely, both  $j(X)$ and $ \bX\smin j(X) $ are unions of the strata.

   Below we will deformation-retract $B$ to its subsets. This deformation-retraction will lift to the deformation-retraction of the strata $\amalg (\bX)_\bullet.$ In particular, this deformation-retraction will be compatible with the splitting $\bX=j(X)\amalg(\bX\smin j(X)) .$

   Therefore below we replace $X$ by $\bX.$ Thus we  assume that the map $X\stackrel{f}{\to}B$ is proper.

\item\ Take a proper morphism $X\to B,$ where  $B$ is a complex affine variety of  $dim_\C B=n$. We construct iterative deformation-retractions of $f.$

Rename  $B$ to $B_n.$  Take the discriminant of $f,$ it is  a Zariski-closed subset   $B_{n-1}\sset B_n,$ of $dim_\C B_{n-1}\le n-1.$
 (This $B_{n-1}$ is   an affine subvariety.)
    The restricted map  $X\smin f^{-1}(B_{n-1})\stackrel{f_|}{\to} B_n\smin B_{n-1}$ is   a topologically locally trivial fibration, by \S\ref{Sec.General.Facts}.iv.

By the Andreotti-Frankel-Hamm theorem $B_n$ is of homotopy type of $dim_\R\le n.$ Therefore, by \S\ref{Sec.General.Facts}.ii, the space $B_n\cup Cells_{n+1}$ is contractible, where $Cells_{n+1}$ is a   union of open cells of $dim_\R \le n+1.$ The gluing of $Cells_{n+1}$ to $B_n$ goes along  the boundaries of these cells. The cells are not algebraic. We can assume that the space $Cells_{n+1}$ is connected and therefore contractible.

 Similarly $B_{n-1}\cup Cells_n$ is contractible, where $Cells_n\sset B_n\cup Cells_{n+1}$ is a union of cells of $dim_\R\le n.$

As $B_{n-1}\sset B_n,$ we have $ B_{n-1}\cap Cells_{n+1}\sseteq \di Cells_{n+1}.$ Thus $dim_\R( B_{n-1}\cap Cells_{n+1})\le n.$
  Moreover, we can assume $dim_\R( B_{n-1}\cap Cells_{n+1})< n,$ as $B_{n-1}$ has no cycles of $dim_\R=n.$

Now we apply \S\ref{Sec.General.Facts}.iii to  the contractible spaces
\beq
Y:=B_n\cup Cells_{n+1},\hspace{1cm} Z:=Cells_{n+1},\hspace{1cm} W:=B_{n-1}\cup Cells_n.
\eeq
Thus $B_n\cup Cells_{n+1}$ deformation-retracts onto $B_{n-1}\cup Cells_n\cup Cells_{n+1}\cup \De_n,$ where $dim_\R\De_n<n+1.$

This deformation-retraction is identity on  $B_{n-1}\cup Cells_n,$ therefore it lifts to a deformation-retraction of the morphism $f.$

\

Finally we check the dimensions:
\bei
\item
If $dim_\R B_{n-1}\le n,$ then $dim_\R (B_{n-1}\cup \De_n\cup Cells_n)\le n.$ Hence the theorem is proved.
\item
Suppose $dim_\R B_{n-1}> n,$ then $dim_\R (B_{n-1}\cup \De_n\cup Cells_n)=dim_\R B_{n-1}.$    Note that $B_{n-1}\cup Cells_n\cup Cells_{n+1}\cup \De_n$ is still contractible. And now iterate the proof.

As the sets $\De_n,$ $Cells_n$ are of $dim_\R\le n,$ we ignore them.
 We consider the map $f^{-1}(B_{n-1})\to B_{n-1}.$
It is proper and algebraic. Take its discriminant, $B_{n-2}\sset B_{n-1}.$ Thus $f$ is a $C^0$-locally trivial fibration over $B_{n-1}\smin B_{n-2}.$   And so on.
\epr
\eei
\eee

\end{document}